\documentclass [12pt] {article}
\title {There are non homotopic framed homotopies of long knots}
\author{Thomas Fiedler}
\usepackage{amssymb,amsfonts,epsfig}
\begin{document}
\newtheorem{proposition}{Proposition}
\newtheorem{theorem}{Theorem}
\newtheorem{lemma}{Lemma}
\newtheorem{corollary}{Corollary}
\newtheorem{example}{Example}
\newtheorem{remark}{Remark}
\newtheorem{definition}{Definition}
\newtheorem{question}{Question}
\newtheorem{conjecture}{Conjecture}
\maketitle
\begin{abstract}
Let $\mathcal {M}$ be the space of all, including singular, long knots in 3-space and for which a fixed projection into the 
plane is an immersion.  Let $cl(\Sigma^{(1)}_{iness})$ be the closure of the union of all singular 
knots in $\mathcal {M}$ with exactly one ordinary double point and such that the two resolutions represent the same (non singular) knot 
type. We call $\Sigma^{(1)}_{iness}$ the {\em inessential walls} and we call 
$\mathcal {M}_{ess} = \mathcal {M} \setminus cl(\Sigma^{(1)}_{iness})$ the {\em essential diagram space}.

We construct a non trivial class in $H^1(\mathcal {M}_{ess}; \mathbb{Z}[A ,A^{-1}])$ by an extension of the 
Kauffman bracket. This implies in particular that there are loops in $\mathcal {M}_{ess}$ which consist of regular 
isotopies of knots together with crossing changings and which are not contractible in $\mathcal {M}_{ess}$ (leading to 
the title of the paper).

We conjecture that our construction gives rise to a new knot polynomial for knots of unknotting number one.

\footnote{2000 {\em Mathematics Subject Classification\/}: 57M25. {\em Key words and phrases\/}:
classical and singular links, discriminant,  knot polynomials, state models}
\end{abstract}

\section{Introduction and results}

The study of knot spaces consists mainly of the study of spaces of non singular knots (compare \cite{BC}, \cite{Bu}, 
\cite{H}, \cite{BCSS}). In this paper we change the point of view: we construct a 1-cocycle for some space which includes
singular knots too. This 1-cocycle is identical zero on all loops which consist only of non-singular knots.

The space of all (possibly singular long) knots  was introduced and studied by Vassiliev in
the pioniering work \cite{V}. It is the space of all differentable maps $f : \mathbb{R} \rightarrow \mathbb{R}^{3}$ which 
agree with $x \rightarrow (x,0,0)$ outside of $[-1,1]$. This is a contractible space (compare \cite{V}), let us call it 
$\mathcal{F}$.
We identify then as usual a knot with its image $f(\mathbb{R})$.
A knot is {\em non-singular} if $f(\mathbb{R})$ is a smooth submanifold, otherwise it is called {\em singular}.
Knots are oriented from the left to the right. 

The singular knots form the discriminant $\Sigma$ of the space of all (long) knots. It has a natural stratification. The complement
of the discriminant are the {\em chambers}. They correspond to the (non-singular) knot types.
Vassiliev has introduced a filtration on the cohomology of the chambers. The 0-dimensional part are the well-known
Vassiliev knot invariants (for all this compare \cite{BN}, \cite{V}, \cite{V2}).

Instead of Vassiliev's knot space $\mathcal{F}$ we will study the space $\mathcal {M}$ of long framed knots.  We fixe a plane in 
3-space which contains the $x-axes$ and we fixe an orthogonal projection $pr$ of the 3-space into
this plane. $\mathcal {M}$ is the space of all those (possibly singular) knots for which the restriction of $pr$ on the knots
is an immersion. The chambers of $\mathcal {M}$ consist of knot diagrams up to isotopy and up to Reidemeister
moves of type II and III (this is called usually {\em regular isotopy} of non singular knots). The points in $\mathcal {M}$
 project to immersed long planar curves. It is easy to see that the connected components of $\mathcal {M}$ are in 1-1 correspondence
 with the regular homotopy classes of immersed long planar curves (standard at infinity). It is well known that the latter 
are in 1-1 correspondence with their Whitney index $n$ (i.e. the degree of the Gauss map). One easily sees that 
$\mathcal {M}$ is a disjoint  union of contractible spaces (numbered by the Whitney index).

The strata of codimension one $\Sigma^{(1)}$ of the discriminant are called the {\em walls}. They correspond to the knots 
which have exactly one ordinary double point as the only singularity. An ordinary double point of an oriented knot 
can be resolved in two different ways (compare Fig. 1).
%Fig1
\begin{figure}
\centering 
\psfig{file=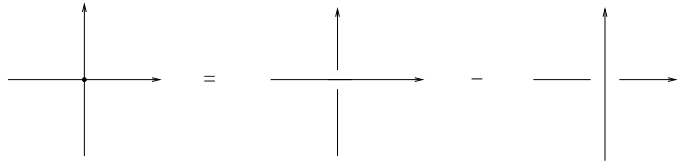}
\caption{Vassiliev's skein relation}
\end{figure}

\begin{definition}
A wall is called {\em inessential} if the two adjacent chambers coincide. The union of all inessential walls is denoted by
$\Sigma^{(1)}_{iness}$. The space  $\mathcal {M} \setminus cl(\Sigma^{(1)}_{iness})$ is
called the {\em essential diagram space} and denoted by $\mathcal {M}_{ess}$. A generic path in $\mathcal {M}_{ess}$
is called a {\em framed homotopy}.
\end{definition}

We could associate to each connected component of $\mathcal {M}$ a stratified space in the following way:
the chambers of $\mathcal {M}$ correspond to vertices. Each stratum $S^{(1)}$ in $\Sigma^{(1)}$ corresponds to an edge 
connecting the vertices corresponding to chambers adjacent to $S^{(1)}$. Each stratum $S^{(2)}$ in $\Sigma^{(2)}$
gives rise to a 2-cell glued to all edges and vertices corresponding to those strata $S^{(1)}$ and those chambers which are 
adjacent to $S^{(2)}$. And we can go on by gluing cells of dimension $n$ corresponding to strata of codimension $n$. 

The resulting space is not a CW-complex because it is not locally finite. (It is well known 
that there are e.g. infinitely many knot types of unknotting number one.) However, this space is evidently still simply 
connected.

Let us now consider the closure $cl(\Sigma^{(1)}_{iness})$ in a component of $\mathcal {M}$ of the union of all 
inessential walls (i.e. we add all adjacent strata of higher codimension). 
 This space is a rather mysterious object. For example, it seems not to be known wether or not
the strata of $\Sigma^{(1)}_{iness}$ correspond always to "nugatory crossings"  (for the definition see e.g. \cite{K1}).
Moreover, I do not know wether or not the complement of $cl(\Sigma^{(1)}_{iness})$ is still connected in each component 
of  $\mathcal {M}$.

\begin{conjecture}
The two knots shown in Fig.2 are not framed homotopic in $\mathcal {M}_{ess}$.
\end{conjecture}.

%Fig2
\begin{figure}
\centering 
\psfig{file=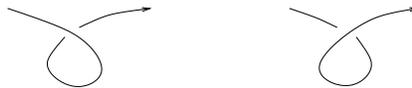}
\caption{not framed homotopic knots?}
\end{figure}

But notice that the two knots shown in Fig.3 are framed homotopic in $\mathcal {M}_{ess}$. The path is shown in 
the figure too. It uses the fact that the "figure eight" knot can be unknotted both by a positive or a negative crossing 
change.

%Fig3
\begin{figure}
\centering 
\psfig{file=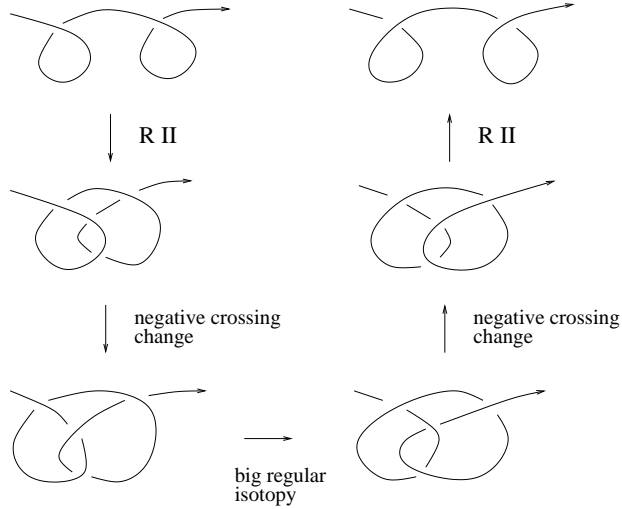}
\caption{a framed homotopy}
\end{figure}

So, we do not know much about the components of the essential diagram space besides the fact that each component
 becomes simply connected (even contractible) if we add again the closure of the non-essential walls 
$cl(\Sigma^{(1)}_{iness})$.

In this paper we describe a surprising phenomen: {\em there are components of the essential diagram space which 
have non-trivial first homology groups}. This result is far from beeing 
obvious and we do not know any other method to prove this.  (In finite dimensions this would imply by Alexander 
duality that $cl(\Sigma^{(1)}_{iness})$ contains a cycle of codimension two in 
$\mathcal {M}$ and this cycle is not homologically trivial in $cl(\Sigma^{(1)}_{iness})$.)

Our prove uses a 1-cocycle which is constructed in the following way:
Let $K$ be a singular knot and let $V^s_K(A,B,C) \in \mathbb{Z}[A ,A^{-1},B,C]$ be the 
extension of the Jones polynomial for singular links contained in \cite{F3} (abusing notation we denote a knot diagram
for $K$ by $K$ too).
It is defined as follows:

$V^s_K(A,B,C) = (-A)^{-3w(K)}<K>_s \in \mathbb{Z}[A ,A^{-1},B,C]$

 (In \cite{F3} we have chosen $C =B^{-1}$ because
the polynomial is homogenous in $B$ and $C$.) The {\em singular Kauffman bracket} $<K>_s$ is defined at crossings 
as the usual Kauffman bracket (compare \cite{K1}) and at double points it is defined as shown in Fig. 4. Here, $B$ and 
$C$ are new independent 
variables. Notice that one of the smoothings induces only a piecewise orientation on the link diagram. $w(K)$ is the writhe
of the knot diagram $K$ (see e.g. \cite{BZ}).
%Fig4
\begin{figure}
\centering 
\psfig{file=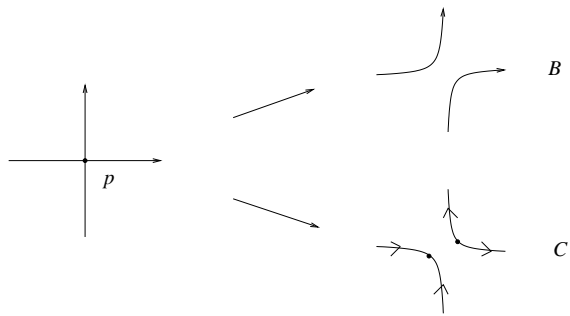}
\caption{smoothings of a double point}
\end{figure}

 Let $\gamma$ be a generic path in $\mathcal{F}$ which connects a given non singular knot with the 
unknot. Here, the end points of the paths are allowed to move inside the chambers. A generic homotopy of such a path 
meets the codimension two part $\Sigma^{(2)}$ of the discriminant
in a finite number of points which correspond to knots with exactly two ordinary double points or with exactly one
ordinary cusp or with exactly one double point with equal tangent directions (see e.g. \cite{C} and \cite{FK}).

 We define a 
{\em co-orientation} on $\Sigma^{(1)}$ by saying that the positive normal direction corresponds to changing a negative
crossing to a positive one. The path $\gamma$ intersects $\Sigma^{(1)}$ transversally in a finite number 
of points.
Let $p$ be such an intersection point. Abusing notation we denote the corresponding double point by $p$ too and
we denote the corresponding singular knot by $K(p)$. Let $ind(p)$ be the intersection index of $\gamma$ with
$\Sigma^{(1)}$ at $p$.

\begin{definition}
Let $\gamma$ be an oriented  path in $\mathcal {M}_{ess}$ which connects a  non singular knot $K$ with  
a diagram of the unknot. The polynomial $Cross(\gamma) \in \mathbb{Z}[A ,A^{-1}]$ is defined by the
following formula

$Cross(\gamma) =  \sum_{p \in \gamma \cap \Sigma} ind(p) <K(p)>_s(A,1, -1)$.
\end{definition}

It follows immediately from the definitions that the intersection index of $\gamma \cap \Sigma$ is an invariant of 
paths $\gamma$ in $\mathcal {M}_{ess}$ up to homotopy in $\mathcal {M}_{ess}$, if 
the end points of the paths are fixed up to regular isotopy. Consequently, the writhe of the unknot at the end of the
path $\gamma$ is completely determined by the writhe $w(K)$ of the knot $K$ and by the intersection index
$ind(\gamma) = \gamma \cap \Sigma$. The Whitney index is invariant under  crossing changes. (As well known, non 
singular knots are regularly isotopic if and only if they are isotopic and they share the same writhe and the same 
Whitney index, compare e.g. \cite{F2}.) In particular, the intersection index of each loop with $\Sigma^{(1)}$ is zero.

\begin{theorem}
Let $\gamma$ be an oriented  loop in $\mathcal {M}_{ess}$. Then the value of $Cross(\gamma)$ depends only on the 
homology class of $\gamma$. The induced cohomology class 
$[Cross] \in H^1(\mathcal {M}_{ess}; \mathbb{Z}[A ,A^{-1}])$ is non trivial.
\end{theorem}

We describe now the first loop $\gamma$ in $\mathcal {M}_{ess}$ for which $Cross(\gamma)$ is non-trivial (we could not 
find a simpler example). 
Let us consider the knots $4_1$ and $6_3$ (from the Knot Atlas) with their 
standard diagrams. They are both amphicheiral. Each of them can be unknotted by either a positive or a negative 
crossing change. Let $K$ be the connected sum  $4_1\#6_3$. Let $h$ be the homotopy which unknots $K$ by a 
positive crossing change (i.e. $ind=1$) of $4_1$ and by a negative crossing change (i.e. $ind=-1$) of $6_3$.
Let $h'$ be the homotopy which unknots $K$ by a negative crossing change of $4_1$ and by a positive
crossing change of $6_3$.

A calculation by hand gives the following result:

$Cross(h\cup(-h')) = (A^{-1}+A+A^{-3}+A^3+A^{-5}+A^5)(<6_3> - <4_1>)$.

Here $<.>$ denotes the usual Kauffman bracket. Notice that the diagrams of $4_1$ and $6_3$ have 
vanishing writhe. Hence, the Kauffman bracket coincides with  the Jones polynomial (see e.g. \cite{K1}, \cite{J}).
We could change $K$ by Reidemeister moves of type I and consider the induced loop $h\cup(-h')$. The value of
$Cross(\gamma)$ 
would change just by a standard factor. This shows that these components of $\mathcal {M}_{ess}$ are non-simply
connected too.

The framed long knot $K$ (i.e. the diagram $K$ or like-wise the smooth submanifold $f(\mathbb{R})$) in our example has 
unknotting number two. The value of  $Cross(\gamma)$ 
looks like the product of some standard 
polynomial (which depends only on the unknotting number of $K$) with a linear combination of Jones polynomials
of those knots of unknotting number one which are contained in the loop. It seems to be very difficult to prove this in general.
But it leads us to the following conjecture.

\begin{conjecture}
Let $K$ be a framed long knot of unknotting number one. Let $h$ be a homotopy which unknots $K$ and such that $ind(h)=1$ 
(respectively $ind(h)=-1$).
Then $Cross(h)$ is invariant under regular  isotopy of $K$.
\end{conjecture}

Notice that apriori the intersection of the closure of the connected component of $K$ in $\mathcal {M}_{ess}$ with the closure 
of the component of the unknot could have an infinite number of connected components. Each of these components 
corresponds to an unknotting of $K$ and they could be all different!

We have verified this conjecture by hand for all diagrams with no more than 6 crossings (which is not a big deal).

\begin{example}

$Cross(3_1, w=3, n=1, ind= -1) = A^{-6} + A^{-4} + A^4$

$Cross(3_1!, w=3, n=1, ind= 1) = -A^{14} - A^{22} - A^{24}$

$Cross(4_1, w= 0, n=0, ind= 1) = A^{-7} + A^{-5} + A + A^3 - A^7$

$Cross(4_1, w= 0, n=0, ind= -1) =A^{-7} - A^{-3} - A^{-1} - A^5 - A^7$.

Here $w$ is the writhe and $n$ is the Whitney index 
of the diagram which we unknot. $3_1$ is the right trefoil and $3_1!$ is its mirror image, the left trefoil.

\end{example}

Our theorem shows that the situation is more complicated already for framed knots of unknotting number two. The corresponding 
conjecture would be that $Cross(h)$ is now well defined up to combinations of Jones polynomials of framed knots with
unknotting number one.

\begin{remark}
There is an extension of the HOMFLY-PT polynomial for singular links by Kauffman and Vogel \cite{KV}. A natural idea 
would be to replace in the construction of $Cross(\gamma)$ the singular Kauffman bracket $<K>_s$ by their singular HOMFLY-PT
polynomial. Surprisingly, this fails. We can not replace $<K>_s$ by our singular Alexander polynomial $\Delta^s$ neither 
(see \cite{F3} and the next section). 
\end{remark}

\begin{remark}
Crossing an inessential wall does not change the knot type but it changes the framing. If we would replace framed 
isotopy by isotopy then our construction would just lead to a version of the Jones polynomial (as it  happened in the first
 version of the present paper).
\end{remark}

\section{Proofs}

Let $\gamma$ be a generic path in the essential diagram space $\mathcal {M}_{ess}$ (i.e. a framed homotopy).

First we observe that $Cross(\gamma)$ is invariant under all homotopies of $\gamma$ which do not change the
intersection with the strata of the discriminant $\Sigma$. This comes from the fact that $<K(p)>_s(A,1, -1)$ is
an isotopy invariant of framed singular knots (i.e. only Reidemeister moves of type II and III and the additional moves
for singular links $SII$ and $SIII$ are allowed, compare \cite{K2} and also \cite{F3}).

Next we need some simple facts from singularity theory which can be proven with the same methods as e.g. in the 
Appendix of \cite{FK} or \cite{C}. 
  Only the following {\em accidents} can occure in a generic homotopy of a path $\gamma$ in $\mathcal {M}_{ess}$
and they can occure only a finite number of times.

I. $\gamma$ becomes tangential (in an ordinary tangent point) to a stratum of $\Sigma^{(1)}$.

II. $\gamma$ passes transversally to a stratum of codimension 2, which consists of the transverse intersection
of two strata of codimension 1. We denote these strata  by $\Sigma^{(2)}_{++}$.

III.$\gamma$ passes transversally to a stratum of codimension 2, which consists of a (single) double point where the 
two branches are tangential. We denote these strata by $\Sigma^{(2)}_{tang}$.

IV. $\gamma$ passes transversally through an ordinary triple point.

All other paths $\gamma$ in the homotopy are just generic paths in $\mathcal {M}_{ess}$.

The important point is that the following accident can occure in $\mathcal {F}$ but not in $\mathcal {M}_{ess}$:

V. $\gamma$ passes transversally to a stratum of codimension 2, which consists of an ordinary cusp.
We denote these strata by $\Sigma^{(2)}_{<}$.

The strata $\Sigma^{(2)}_{<}$ are just (a generic part of) boundaries of inessential walls in $\mathcal {F}$, but all our 
paths are in $\mathcal {M}_{ess}$. 

%The main degenerations of $\Sigma^{(2)}_{++}$ are ordinary triple points and a couple of double points in a self-tangency.
% They form strata
%of codimension 3. Triple points were already studied in order to obtain the well known
%{\em 4T-relation} in the theory of finite type invariants (compare e.g. \cite{BN}). 

%We do not need the strata of $\Sigma^{(3)}$ in this paper.
%But we have decided to add a picture of the meridional 2-sphere of a stratum corresponding to a triple point, because 
%we have not seen the picture in this form before.
%It could become important in connection with our Remark 4. 
%Fig3
%\begin{figure}
%\centering 
%\psfig{file=cross3.eps}
%\caption{unfolding of a triple point}
%\end{figure}

%The three circles on the 2-sphere in Fig. 3 are the intersections with $\Sigma^{(1)}$. They intersect in six strata of 
%$\Sigma^{(2)}_{++}$, which come in antipodal pairs. In Fig. 3 we describe the closure of just one of the eight adjacent
%chambers. This is enough to establish the rest of the picture.
%At a generic triple point the three tangent vectors span a 3-dimensional space. But notice, that if we fix in addition a generic
%projection into a plane  (as usual in knot theory) then we have to distinguish {\em braid-like} triple
%points from {\em star-like} triple points (compare \cite{F2}).

It follows immediately from the definitions that all our polynomials are invariant under accidents of type I.

In order to prove Theorem 1 we have to show that $Cross(m) = 0$ for the meridional loops $m$ of $\Sigma^{(2)}_{++}$
 and of $\Sigma^{(2)}_{tang}$.

%Fig5
\begin{figure}
\centering 
\psfig{file=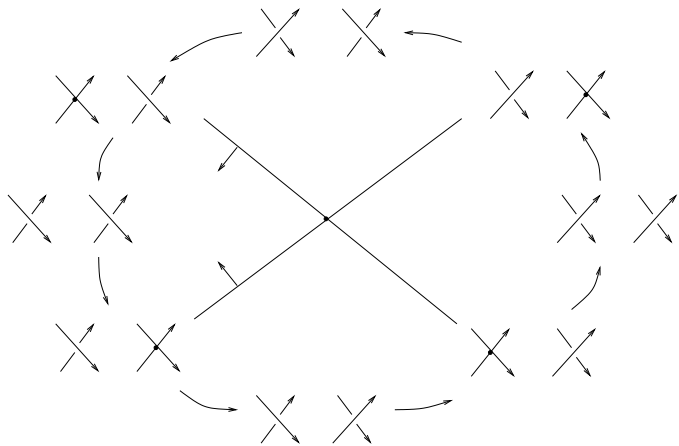}
\caption{meridian for a pair of double points}
\end{figure}

%Fig6
\begin{figure}
\centering 
\psfig{file=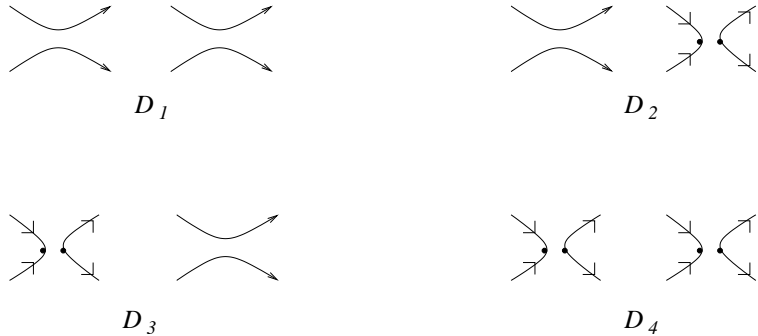}
\caption{smoothed diagrams near the double points}
\end{figure}

We show the diagrams for the meridional loop of $\Sigma^{(2)}_{++}$ in Fig. 5.  After smoothing the crossings and
the double points we are left with exactly the four diagrams $D_i$ shown in Fig. 6. They are in general independent. 
The meridional loop 
$m$ of $\Sigma^{(2)}_{++}$ leads  to the following equations:

$(-AB-A^{-1}B+AB+A^{-1}B)<D_1> = 0$

$(-A^{-1}B-A^{-1}C+AC+AB)<D_2> = 0$

$(-AC-AB+A^{-1}B+A^{-1}C)<D_3> = 0$

$(-A^{-1}C-AC+A^{-1}C+AC)<D_4> = 0$.

Here, $<D_i>$ is the usual Kauffman bracket. $Cross(m) =0$ in general  if and only if each coefficient of $<D_i>$ is 
zero. Therefore we obtain the unique (non trivial) solution

$C= -B$.

Each singular link in the framed homotopy has exactly one double point. Therefore we do not loose information by 
setting $B=1$,
and we obtain exactly the definition of $Cross(\gamma)$.
%Fig7
\begin{figure}
\centering 
\psfig{file=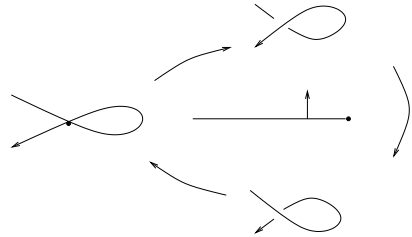}
\caption{meridian for a cusp}
\end{figure}

Notice that it could happen that exactly one of the four strata of $\Sigma^{(1)}$ in Fig. 5 is non essential. In this case
we could not push our path $\gamma$ in $\mathcal {M}_{ess}$ through the corresponding stratum of 
$\Sigma^{(2)}_{++}$ by a small homotopy.
Our theorem implies that sometimes we can not even do it by a big homotopy.

We show the diagrams for the meridional loop of $\Sigma^{(2)}_{<} \subset \mathcal{F}$ in Fig. 7. 
The value of $Cross(m)$ is a non trivial multiple of a Kauffman bracket. 
Consequently, $Cross(m)$ is in general non zero.

%Fig8
\begin{figure}
\centering 
\psfig{file=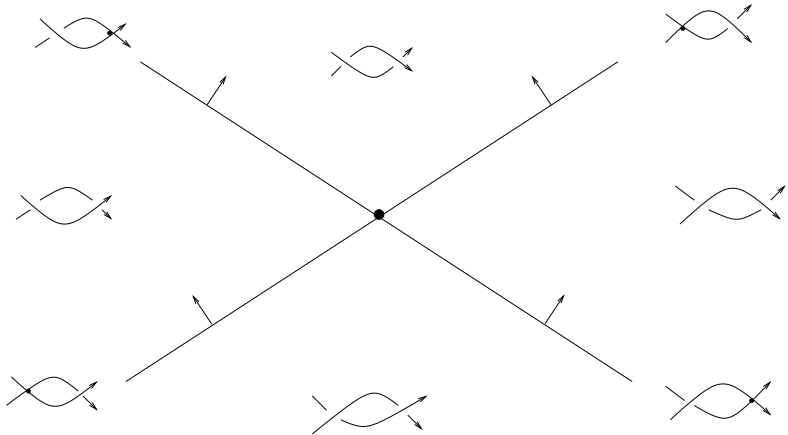}
\caption{meridian for a double point with equal tangencies}
\end{figure}

We show the diagrams for the meridional loop of $\Sigma^{(2)}_{tang}$ in Fig. 8.  

$Cross(m) =0$ follows from the
invariance of $V^s$ under the move $S II$ for singular links (compare \cite{F3}). (We had to choose in Fig. 8 a type
of a Reidemeister II move and orientations on the branches. Taking mirror images or changing orientations of
branches leads to the same result.)

The stratum of an ordinary triple point is not smooth. But one easily sees that the meridians form  theta-graphs.
Using the invariance of $Cross$ under passing $\Sigma^{(2)}_{tang}$ it suffices to prove the invariance  for the loops in  just one of the 
theta graphs (compare \cite{F2}).
The three arcs $m_1,m_2,m_3$ in the theta-graph are shown in Fig. 10,11 and 12. It suffices hence to prove that 
$Cross(m_1 \cup (-m_2)) = Cross(m_1 \cup (-m_3)) = 0$.

%Fig9
\begin{figure}
\centering 
\psfig{file=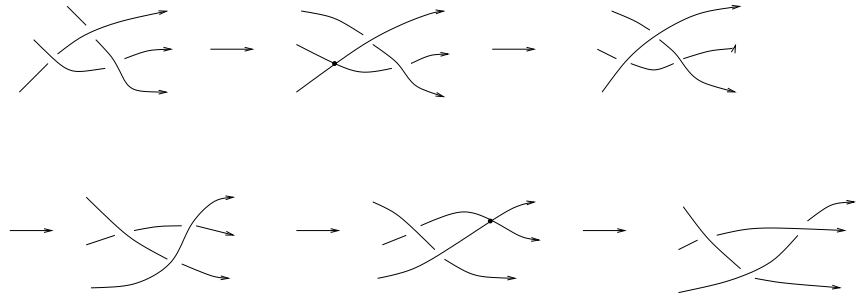}
\caption{the path $m_1$}
\end{figure}
%Fig10
\begin{figure}
\centering 
\psfig{file=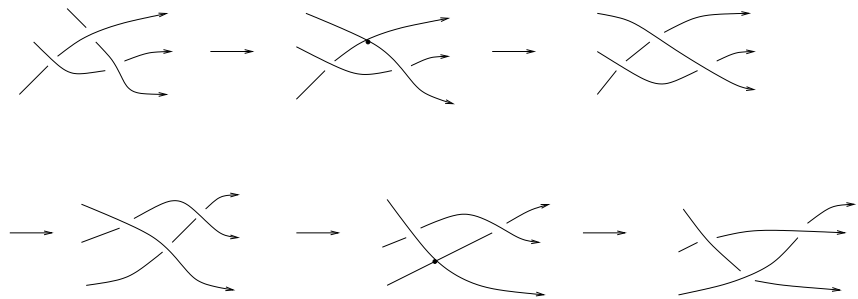}
\caption{the path $m_2$}
\end{figure}
%Fig11
\begin{figure}
\centering 
\psfig{file=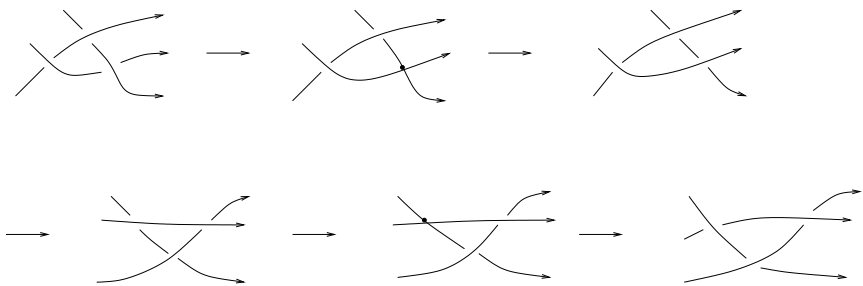}
\caption{the path $m_3$}
\end{figure}
%Fig12
\begin{figure}
\centering 
\psfig{file=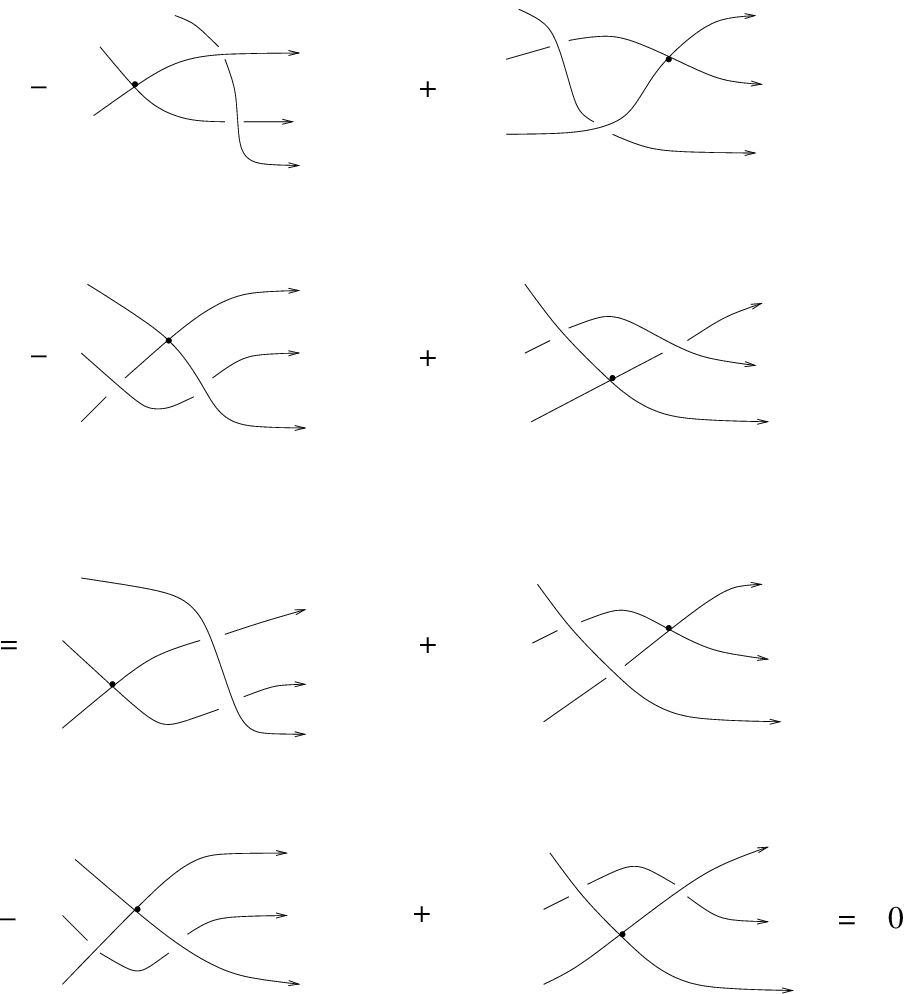}
\caption{$Cross(m_1 \cup (-m_2))$}
\end{figure}

We show the calculation of $Cross(m_1 \cup (-m_2))$ in Fig. 12. The first equality in Fig. 12 uses twice the fact that 
$Cross$ is invariant under passing $\Sigma^{(2)}_{++}$. The calculation of $Cross(m_1 \cup (-m_3))$ is completely
analogous and is left to the reader.

It is clear that we can replace homotopy by homology in all our constructions because the values are in a commutative
ring. 
We have proven that $Cross(\gamma)$ depends only on the homology class of $\gamma$ in $\mathcal {M}_{ess}$.
The example in the introduction shows that $Cross$ induces a non trivial cohomolgy class.
The theorem is proven.
\vspace{0.3cm}

It remains to prove our assertion from Remark 1.
We recall the definition of the singular HOMFLY-PT polynomial from \cite{KV} in Fig. 13. Here, $x, y, a$ are new 
independent variables. We want to replace $V^s$ by this polynomial. A calculation for the meridional loop of 
$\Sigma^{(2)}_{++}$ leads to the equation shown in Fig. 14, which has the unique solution $x = y$. But then the 
singular HOMFLY-PT polynomial is no longer sensitive for crossing changes of links.
%Fig13
\begin{figure}
\centering 
\psfig{file=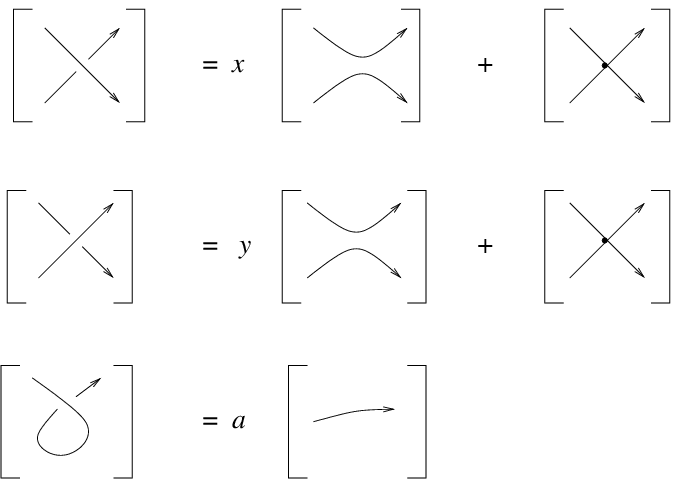}
\caption{Kauffman-Vogel's skein relation}
\end{figure}
%Fig14
\begin{figure}
\centering 
\psfig{file=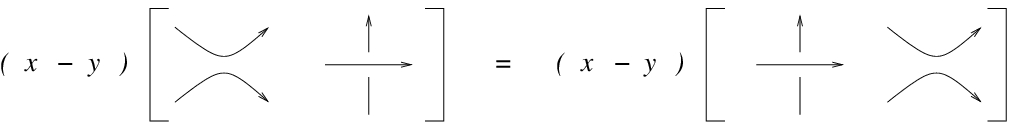}
\caption{equation from the meridian for a pair of double points}
\end{figure}

The failure of the singular Alexander polynomial from \cite{F3} is also caused by the stratum $\Sigma^{(2)}_{++}$.
It leads to the equation $A^4 = 1$, which makes the invariant uninteresting. We leave the verification to the reader.

{\em Acknowledgements}--- I am grateful to the referee for his useful comments.

Institut de Math\'ematiques de Toulouse

Universit\'e Paul Sabatier et CNRS (UMR 5219)

118, route de Narbonne 

31062 Toulouse Cedex 09, France

fiedler@picard.ups-tlse.fr

\end{document}